\documentclass[11pt]{amsart}
\usepackage{palatino}
\usepackage[all]{xy,xypic}
\usepackage{amsfonts,amsmath,oldgerm,amssymb,amscd}

\newcommand{\ra}{\rightarrow}		
\newcommand{\lra}{\longrightarrow}
\newcommand{\by}[1]{\stackrel{#1}{\ra}}

\newcommand{\remove}[1]{}
\newcommand{\e}{{\bf E }}

\newcommand{\surj}{\ra\!\!\!\ra}	

\newcommand{\ol}{\overline}		
\newcommand{\wt}{\widetilde}
\newcommand{\iso}{\by \sim}

\newtheorem{theorem}{Theorem}[section]
\newtheorem{proposition}[theorem]{Proposition}
\newtheorem{lemma}[theorem]{Lemma}
\newtheorem{definition}[theorem]{Definition}
\newtheorem{corollary}[theorem]{Corollary}
\newtheorem{conjecture}[theorem]{Conjecture}
\newtheorem{example}[theorem]{Example}

\newcommand{\ga}{\alpha}	\newcommand{\gb}{\beta}
		\newcommand{\gd}{\delta}

\newcommand{\bq}{\mbox{$\mathbb Q$}}

	\newcommand{\bz}{\mbox{$\mathbb Z$}}

\newcommand{\CI}{\mbox{$\mathcal I$}}

\newcommand{\mm}{\mbox{$\mathfrak m$}}

\newcommand{\ot}{\mbox{\,$\otimes$\,}}	
\newcommand{\op}{\mbox{$\oplus$}}
\newcommand{\Spec}{\text{Spec}}	
\newcommand{\hh}{\text{ht}}

\newcommand{\dd}{\text{dim}}

\newcommand{\sur}{\twoheadrightarrow}
\newcommand{\bp}{\begin{proposition}}
\newcommand{\ep}{\end{proposition}}
\newcommand{\bl}{\begin{lemma}}
\newcommand{\el}{\end{lemma}}
\newcommand{\bt}{\begin{theorem}}
\newcommand{\et}{\end{theorem}}
\newcommand{\bc}{\begin{corollary}}
\newcommand{\ec}{\end{corollary}}
\newcommand{\bd}{\begin{definition}}
\newcommand{\ed}{\end{definition}}
\newcommand{\bco}{\begin{conjecture}}
\newcommand{\eco}{\end{conjecture}}

\renewcommand{\l}{\longrightarrow}

\def\rmk{\refstepcounter{theorem}\paragraph{{\bf Remark} \thetheorem}}
\def\proof{\paragraph{Proof}}

\def\notation{\paragraph{\bf Notation}}
\def\quest{\refstepcounter{theorem}\paragraph{{\bf Question} \thetheorem}}
\def\definition{\refstepcounter{theorem}\paragraph{{\bf Definition} \thetheorem}}
\title[]{Remarks on Euler class groups and two conjectures}
\author{ Mrinal Kanti Das}
\address{Stat-Math Unit, Indian Statistical Institute, 203 B. T. Road, Kolkata 700108 India}
\email{mrinal@isical.ac.in}
\date{\today}
\thanks{}

\begin{document}
\maketitle

%\begin{abstract}
 %\end{abstract}

%\maketitle

\section{Introduction}
Recently,  Asok-Fasel have settled two  fascinating open problems in \cite{af}. Let us  recall their results.
Let $k$ be a field and $X=\Spec(R)$ be a smooth affine $k$-scheme  of dimension $d\geq 2$. Let 
$E^d(R)$ be the $d\,\textsuperscript{th}$ Euler class group of Nori-Bhatwadekar-Sridharan \cite{brs1} and 
$\wt{CH}^d(X)$ be the Chow-Witt group defined by Barge-Morel \cite{bm} (and studied by Fasel \cite{f} in detail). Also, consider the Chow group $CH^d(X)$ of zero cycles, and the weak Euler class group $E_0^d(R)$ of Bhatwadekar-Sridharan (introduced in \cite{brs2}).   Asok-Fasel prove the following comparison theorems.

\bt\label{af}\cite[Theorems 3.2.1,  3.2.7]{af}
Let $k$ be a field and $X=\Spec(R)$ be a smooth affine $k$-scheme  of dimension $d\geq 2$.
Assume that: (a) $k$ has characteristic zero if $d=3$, or (b) $\text{char}(k)\neq 2$, $k$ is infinite and perfect if $d=2$ or 
$d\geq 4$.  Then,
\begin{enumerate}
\item
 there is a functorial isomorphism $E^d(R)\iso \wt{CH}^d(X)$;
\item
the canonical morphism $E_0^d(R)\lra CH^d(X)$ is an isomorphism. 
\end{enumerate}
\et

We shall remark that the above results can  easily be  extended to the following form.

\bt\label{our}
Let $k$ be an infinite perfect field and $X=\Spec(R)$ be a smooth affine $k$-scheme  of dimension $d\geq 2$. Then,
\begin{enumerate}
\item
if $\text{char}(k)\neq 2$, there is a functorial isomorphism $E^d(R)\iso \wt{CH}^d(X)$;
\item
if $d=2$, or $d\geq 3$ and $\text{char}(k)\neq 2$, the canonical morphism $E_0^d(R)\lra CH^d(X)$ is an isomorphism. 
\end{enumerate}
\et

\rmk
The proof of $E_0^d(R)\iso CH^d(X)$ (for $d\geq 3$) relies on the isomorphism $E^d(R)\iso \wt{CH}^d(X)$ and
the assumption $\text{char}(k)\neq 2$ stems from there. Therefore, the following remains open.

\quest
Can one remove the characteristic assumption from Theorem \ref{our} (2) for $d\geq 3$? 

\medskip

The  improvements of the results of Asok-Fasel in Theorem \ref{our}, despite appearing significant,  did not warrant much work. The ingredients are all implicit in the existing theory of the Euler class groups as developed in \cite{brs1,brs2,brs3,d1}. We have only  highlighted them properly and made them explicit. 

Let us explain to the reader why the odd assumption  \emph{``$k$ has characteristic zero if $d=3$"} has landed in the hypothesis of Theorem \ref{af}. In \cite[Lemma 3.1.11]{af} Asok-Fasel have used the \emph{homotopy invariance}
of the Euler class groups ($E^d(R)\iso E^d(R[T])$). For $d\geq 4$ they refer to \cite{my}, while for $d=3$ they have to use \cite{d1}. In \cite{d1}, the theory of the Euler class group $E^d(R[T])$ ($d\geq 3$) was developed with the blanket assumption that $\bq\subset R$, and for this reason Asok-Fasel assume that $\text{char}(k)=0$ in Theorem \ref{af}. To rectify this, here we define $E^d(R[T])$ for $d\geq 3$ where $R$ is a regular ring containing a field  $k$
($R$ is of finite type over $k$ if $k$ is finite). Then we
indicate how the homotopy invariance can be achieved when $R$ is a smooth affine algebra of dimension $d\geq 3$ over an infinite field $k$. In our definition of $E^d(R[T])$,  for an Euler cycle $(I,\omega_I)\in E^d(R[T])$ we take the \emph{local orientation} $\omega_I$ as an equivalence class of surjections induced by the action of
$SL_d(R[T]/I)$ on all surjections from $(R[T]/I)^d\sur I/I^2$. On the other hand, the definition followed by  \cite{my} (from \cite{brs4}) uses the
action of the elementary group $E_d(R[T]/I)$ (which may be a proper subgroup of $SL_d(R[T]/I)$). Therefore, a priori the two definitions are different and it takes some amount of work to show their equivalence (one has to use Theorem \ref{sl} from below).

Our starting point of this article was what we just described above. But
while carrying out the above task, we felt that we should also revamp some parts of the existing theory of Euler classes and present them in as much generality as possible so that others can use them if necessary.  On the way, we have obtained some new  results as well (see \ref{new1}, \ref{new2}, \ref{new3}) which will perhaps attract attention of a general reader and will find some interesting applications. Again, it will be apparent to an expert that our methods involve hardly any substantial original idea.

\section{The Euler class group $E^d(R)$}\label{ecg}
The Euler class group $E^d(R)$ can be defined for any commutative Noetherian ring $R$ of dimension $d\geq 2$.  Let us first recall the two definitions of $E^d(R)$ from \cite{brs1} and \cite{brs3}.
But before doing so, we insert a definition which will be frequently used.

\smallskip

\definition
({\bf Local and global orientations:})
Let $R$ be a commutative Noetherian ring  of 
dimension $d\geq 2$.
Let $J\subset R$ be an ideal of height $d$ such that $J/J^2$ is generated 
by $d$ elements. Two surjections $\alpha,\beta$ from $(R/J)^d$ to $J/J^2$ are
said to be related if there exists $\sigma \in SL_{d}(R/J)$ such that 
$\alpha\sigma=\beta$. Clearly this is an equivalence relation on the set
of surjections from $(R/J)^d$ to $J/J^2$. Note that, if $\alpha$ can be lifted to a 
surjection $\theta :R^d\sur J$, then $\beta$ can also be lifted to some
$\theta':R^d\sur J$.
Let $[\alpha]$ denote the equivalence 
class of $\alpha$. Such an equivalence class $[\alpha]$
is called a \emph{local orientation} of $J$. By abuse of notation, 
we shall identify an equivalence class $[\alpha]$ with $\alpha$.
A local orientation $\alpha$ is called a \emph{global orientation}
if $\alpha :(R/J)^d\sur J/J^2$ can be lifted to a surjection $\theta :R^d\sur J$.

\smallskip

\definition\label{def1}
({\bf The Euler class group  $E^d(R)$: $R$ smooth affine}) Let $R$ be a smooth affine domain of dimension $d\geq 2$
over an infinite perfect field $k$. Let $G_1$ be the free abelian group on the set of pairs $(\mm,\omega_{\mm})$, 
where $\mm$ is a maximal ideal of $R$ and $\omega_{\mm}$ is a local orientation of $\mm$. Let $J=m_{1}\cap \cdots \cap m_r$, where $m_i$ are distinct  maximal
ideals of $R$. Any $\omega_{J}:(R/J)^d\surj J/J^2$ induces surjections
$\omega_i :(R/m_i)^d\surj m_i/m_i^2$ for each $i$. We associate
$(J,\omega_J):= \sum_{1}^{r}(m_i,\omega_i)\in G_1$. Let $S_1$ be
the set of elements $(J,\omega_J)$ of $G_1$ for
which $\omega_J$ is a global orientation and $H_1$ be the 
subgroup of $G_1$ generated by $S_1$ .  The
Euler class group $E^d(R)$ is defined as $E^d(R))\stackrel{\text{ def}}{=}G_1/H_1$.

\smallskip

\definition\label{def2}
({\bf The Euler class group  $\wt{E}^d(R)$: $R$ Noetherian})  Let $R$ be a commutative Noetherian ring of dimension $d\geq 2$.
Let $G_2$ be the free abelian group on the set of pairs 
$({\mathcal N},\omega_{\mathcal N})$ where ${\mathcal N}$ is an 
${\mathcal M}$-primary ideal for some maximal ideal ${\mathcal M}$ of height
$d$ such that  ${\mathcal N}/{\mathcal N}^2$ is generated by $d$ elements
and $\omega_{\mathcal N}$ is a local orientation of ${\mathcal N}$.
Now let $J\subset R$ be an ideal of height $d$ such that $J/J^2$ is 
generated by $d$ elements
and $\omega_{J}$ be a local orientation of $J$. Let 
$J=\cap_{i}{\mathcal N}_{i}$ be  the (irredundant) primary decomposition
of $J$. We associate to the pair $(J,\omega_{J})$, the element 
$\sum_{i} ({\mathcal N}_{i},\omega_{{\mathcal N}_{i}})$ of $G_2$ where 
$\omega_{{\mathcal N}_{i}}$ is the local orientation of ${\mathcal N}_{i}$
induced by $\omega_{J}$. By abuse of notation, we denote 
$\sum_{i} ({\mathcal N}_{i},\omega_{{\mathcal N}_{i}})$
by $(J,\omega_{J})$.
Let $H_2$ be the subgroup of $G_2$ generated by the set of pairs $(J,\omega_{J})$,
where $J$ is an ideal of height $d$  and $\omega_{J}$ is a global
orientation of $J$.
The Euler class group of $R$  is 
$\wt{E}^d(R)\stackrel{\text{ def}}{=}G_2/H_2$.   

To explain that the above two notions coincide when $R$ is a smooth affine domain of dimension $d\geq 2$
over an infinite perfect field $k$, we have the following proposition.

\bp
Let $R$ be a smooth affine domain of dimension $d\geq 2$
over an infinite perfect field $k$. Then $E^d(R)\simeq \wt{E}^d(R)$.
\ep

\proof
The map which sends each $(\mm,\omega_{\mm})$ to itself gives rise to a morphism from $G_1$ to $G_2$.
Since this morphism takes $H_1$ inside $H_2$, we get a canonical morphism $\phi:E^d(R)\lra \wt{E}^d(R)$.
As remarked in \cite[Remark 4.16]{brs1}, an element of $E^d(R)$ is represented by a single Euler cycle $(J,\omega_J)$,
where $J$ is a reduced ideal of height $d$ and $\omega_J$ is a local orientation. If $\phi((J,\omega_J))=
(J,\omega_J)=0$ in $\wt{E}^d(R)$, then by \cite[Theorem 4.2]{brs3} it follows that $\omega_J$ is a global orientation and therefore, $(J,\omega_J)=0$ in $E^d(R)$. Therefore, $\phi$ is injective. Now let $(I,\omega_I)$ be an element of $\wt{E}^d(R)$. As $k$ is infinite and perfect, we can apply Swan's Bertini theorem \cite[Theorem 2.11]{brs2} to find $(I',\omega_{I'})\in \wt{E}^d(R)$ such that $I+I'=R$, $I'$ is reduced of height $d$,  and
$(I,\omega_I)+(I',\omega_{I'})=0$ in $\wt{E}^d(R)$ (see \cite[2.7, 2.8]{drs1} for a proof). As $(I',\omega_{I'})$ is in the image of $\phi$, it follows that 
$\phi$ is surjective.
\qed

\smallskip

From now on, both the groups will be denoted by $E^d(R)$.

 \section{The Euler class of a projective module}\label{proj}
 Let $R$ be a commutative Noetherian ring of dimension $d\geq 2$.
Let $P$ be a projective $R$-module of rank $d$ such that
$ R \simeq \wedge^d (P)$ and let 
$\chi : R\iso \wedge^{d}P$ be an isomorphism. Let 
$\varphi: P \sur J$ be a surjection where $J$ is an ideal of height $d$.
Therefore we obtain an induced surjection 
${\overline \varphi}:  P/JP \sur J/J^2.$ As $\dd(R/J)=0$, we see that $P/JP$ is a free
$R/J$-module of rank $d$.
Let $\overline \gamma: (R/J)^{d}
\simeq P/JP$ be an isomorphism  such  that $\wedge^d
(\overline \gamma ) = \overline \chi .$ Let $\omega_J $ be the 
local orientation of $J $ given by  $ \overline
{\varphi }~ \overline \gamma : (R/J)^{d} \sur
 J/J^2.$ Let  $e_d(P,\chi)$ be the image in $E^d(R)$
of the element $(J, \omega_J )$. 
The \emph{Euler class} of $(P,\chi)$  is defined to be $e_d(P,\chi)$,
\emph{provided, the assignment sending the pair $(P, \chi) $ to the element
$e_d(P,\chi)$ of  $E^d(R)$  is independent of the choice of the surjection $\varphi: P \sur J$}.

\medskip

\quest
Is $e_d(P,\chi)$, as described above, well-defined?

\medskip

We do not know the complete answer yet, when $R$ is just commutative Noetherian.
Before showing  that $e_d(P,\chi)$ is well-defined for certain classes of rings, 
we first make the following definition.

\smallskip

\definition\label{prope}
A ring $R$ is said to have property {\bf E} if every projective $R[T]$-module is extended from $R$.

In this context, we  may recall the following conjecture.

\smallskip

\begin{conjecture}\label{bqc}(Bass-Quillen) 
Any regular ring has property {\bf E}.
\end{conjecture}

\medskip

We shall comment more on this later. Let us also recall a very useful result by Roitman
\cite[Proposition 2]{ro}.

\bt\label{roit}
Let $R$ be a ring and $S$ be a  multiplicative subset of $R$. If $R$ has property {\bf E} then so does $R_S$.
\et

We shall need the following \emph{``moving lemma"}. The version given below can easily be proved following \cite[2.4]{brs4}, which in turn is essentially  based on \cite[2.14]{brs3}.

\bl\label{move}
Let $A$ be a Noetherian ring of dimension $d$ and let $J\subset A$ be an ideal of height $n$
such that $2n\geq d+1$. Let $P$ be a projective $A$-module of rank $n$ and $\ol{\ga}:P/JP\sur J/J^2$ be a surjection. Then, there exists an ideal $J'$ of $A$ and a surjection $\gb:P\sur J\cap J'$ such that:
\begin{enumerate}
\item
$J+J'=A$,
\item
$\gb \ot A/J=\ol{\ga}$,
\item
$\hh(J')\geq n$.
\end{enumerate}
Given any ideal $K\subset A$ of height $n$, the map $\gb$ can be  chosen so that $J'+K=A$.  
\el

\rmk
In the above lemma, if $n=d$ and if $A$ is a geometrically reduced affine algebra over an infinite field $k$, then using Swan's Bertini Theorem one can ensure that the ideal $J'$ (if proper) is reduced.

\medskip

The proof of the following theorem is based on that  of \cite[5.2]{brs4}.

\bt\label{chi}
Let $A$ be a Noetherian ring of dimension $d\geq 2$ with the following  assumptions:
\begin{enumerate}
\item[(a)]
if $d=2$: no restriction;
\item[(b)]
if $d=3$: either $A$ is a $\bz[\frac{1}{2}]$-algebra, or $A$ is  regular; 
\item[(c)]
$d\geq 4$: either $A$ is a $\bz[\frac{1}{(d-1)!}]$-algebra, or $A$ has property {\bf E}.
\end{enumerate} 
Let $P$ be a projective $A$-module of rank $d$ and let  $\ga(T):P[T]\sur I$ be a surjection where $I$ is an ideal of $A[T]$ height $d$. Assume that 
$J=I(0)$ is a proper ideal, and further that $P/NP$ is free,  where $N=(I\cap A)^2$. Let $p_1,\cdots,p_{d}\in P$ be such that their images in $P/NP$ form a basis. Let $a_1,\cdots,a_{d}\in J$ be such that 
$\ga(0)(p_i)=a_i$. Then, there exists an ideal $K\subset A$ of height $\geq d$ such that $K+N=A$ and:

\begin{enumerate}
\item
$I\cap K[T]=(F_1(T),\cdots,F_{d}(T))$,
\item
$F_i(0)-F_i(1)\in K^2$, $i=1,\cdots,d$,
\item
$\ga(T)(p_i)-F_i(T)\in I^2$, $i=1,\cdots,d$,
\item
$F_i(0)-a_i\in J^2$, $i=1,\cdots,d$.
\end{enumerate}
\et

\proof
The proof follows \cite{brs4} up to a certain part. Then,
modifications are only needed in \cite[page 151, second paragraph onward]{brs4}. However, we write the whole proof for completeness.

Since $\ga(0)(P/NP)=J/NJ$, it follows that $J=(a_1,\cdots,a_d)+NJ$. Since $NJ\subset J^2$, there exists
$c\in NJ$ such that $J=(a_1,\cdots,a_d,c)$. Using a standard general position argument and Lemma \ref{move}, we can find $b_1,\cdots,b_d\in A$ such that:
\begin{enumerate}
\item[(a)]
$(a_1+cb_1,\cdots,a_d+cb_d)=J\cap K$, where $\hh(K)\geq d$;
\item[(b)]
$K+(c)=A$ and therefore, $K+N=A=K+J$. 
\end{enumerate}

As $c\in NJ$, we have, $c=\sum c_jd_j$, where $c_j\in N$ and $d_j\in J$. Let $q_j\in P$ be 
such that $\ga(0)(q_j)=d_j$ (the map $\ga(0):P\sur J$ is surjective). Let $\wt{p_i}=p_i+b_i\sum c_jq_j$.
Then $\wt{p_i}-p_i\in NP$. Let $\ga(T)(\wt{p_i})=G_i(T)\in I$. We record that $G_i(0)=a_i+cb_i$.

Let $I'=I\cap K[T]$. Then $I'/I'^2=I/I^2\op K[T]/K^2[T]$ (as $I+K[T]=A[T]$), and $I'(0)=J\cap K=(a_1+cb_1,\cdots,a_d+cb_d)$.

As $P/NP$ is free, it follows that $P_{1+N}$ is a free $A_{1+N}$-module with basis $\wt{p_1},\cdots,\wt{p_d}$. We have,
$$I'_{1+N}=I_{1+N}=(G_1(T),\cdots,G_d(T)).$$
We can choose $a\in N$ such that $1+a\in K$ and $I'_{1+a}
=I_{1+a}=(G_1(T),\cdots,
G_{d}(T))$. On the other hand, $I'_a=K_a[T]=(a_1+cb_1,\cdots,a_{d}+cb_{d})$ such that 
$G_i(0)=a_i+cb_i$ for $i=1,\cdots,d$.

We now split the cases.

\smallskip

\noindent
{\it Case 1.} Assume that $d=2$. Let $b=1+a$. The rows $ (G_1(T),G_{2}(T))$ and $(a_1+cb_1,a_{2}+cb_{2})$ are unimodular over the ring 
$A_{ab}[T]$, and they agree when $T$ is set to zero.
As any unimodular row of length two over any ring can be completed to a $2\times 2$ with determinant one, we can find $\theta(T)\in SL_2(A_{ab}[T])$ such that $(G_1(T),G_2(T))\theta(T)=(1,0)$. Then $(G_1(0),G_2(0))\theta(0)=(1,0)$, implying
that $(a_1+cb_1,a_2+cb_2)\theta(0)=(1,0)$. Taking $\sigma(T):=\theta(T)\theta(0)^{-1}$, we
observe that $\sigma(0)=id$, and $(G_1(T),G_2(T))\sigma(T)=(a_1+cb_1,a_2+cb_2)=(G_1(0),G_2(0))$. The rest of the arguments are the same as \cite[5.2]{brs4}.

\smallskip

\noindent
{\it Case 2.} Assume that $d=3$.   Consider the two unimodular rows $(G_1(T),G_2(T),G_{3}(T))$
and $(a_1+cb_1,a_2+cb_2,a_{3}+cb_{3})$ over the ring $A_{a(1+aA)}[T]$. Note that $\dd(A_{a(1+aA)})\leq 2$. 
If $A$ is a $\bz[\frac{1}{2}]$-algebra, then, by a result of Murthy \cite[Theorem 2.5]{ra1}, the unimodular row $(G_1(T),G_2(T),G_{3}(T))$ is locally completable, and therefore by Quillen's local-global principle, it is extended from $A_{a(1+aA)}$. On the other hand, if $A$ is a regular ring, then 
by  a result of Murthy \cite{mu1}, the unimodular row $(G_1(T),G_2(T),G_{3}(T))$ is locally completable, and therefore by Quillen's local-global principle, it is extended from $A_{a(1+aA)}$. In other words, there exists $\theta(T)\in GL_3(A_{a(1+aA)}[T])$ such that 
$$(G_1(T),G_2(T),G_{3}(T))\theta(T)=(G_1(0),G_2(0),G_3(0))=(a_1+cb_1,a_2+cb_2,a_{3}+cb_{3})$$
Taking $\sigma'(T):=\theta(T)\theta(0)^{-1}$, we
observe that $\sigma'(T)\in SL_3(A_{a(1+aA)}[T])$, $\sigma'(0)=id$, and $(G_1(T),G_2(T),G_{3}(T))\sigma'(T)=(a_1+cb_1,a_2+cb_2,a_{3}+cb_{3})$. 
 We can find  some $b$ of the form $1+\lambda a$ such that $b$ is a multiple of $1+a$, and 
some $\sigma(T)\in SL_3(A_{ab}[T])$ such that $\sigma(0)=id$ and $(G_1(T),G_2(T),G_{3}(T))\sigma(T)=
(a_1+cb_1,a_2+cb_2,a_{3}+cb_{3})$ over the ring $A_{ab}[T]$. The rest is same
as \cite[5.2]{brs4}.

\smallskip

\noindent
{\it Case 3.} Assume that $d\geq 4$.  Consider the unimodular rows $(G_1(T),\cdots,G_{d}(T))$
and $(a_1+cb_1,\cdots,a_{d}+cb_{d})$ over the ring $A_{a(1+aA)}[T]$. They agree on $T=0$. Now, if $A$ is a $\bz[\frac{1}{(d-1)!}]$-algebra, then by \cite[Theorem 2.4]{ra2},the row $(G_1(T),\cdots,G_{d}(T))$
is locally completable,  and therefore by Quillen's local-global principle, it is extended from 
$A_{a(1+aA)}$. On the other hand,
if $A$ has property {\bf E}, then by Theorem \ref{roit}, so does $A_{a(1+aA)}$. Therefore, $(G_1(T),\cdots,G_{d}(T))$
is extended from $A_{a(1+aA)}$.
Same line of arguments as in Case 2 will work.
\qed

\medskip

We now prove:

\bt\label{ecp}
Let $R$ be a Noetherian ring of dimension $d\geq 2$ with the following  assumptions:
\begin{enumerate}
\item[(a)]
if $d=2$: no restriction;
\item[(b)]
if $d=3$: either $R$ is a $\bz[\frac{1}{2}]$-algebra, or $R$ is  regular;
\item[(c)]
$d\geq 4$: either $R$ is a $\bz[\frac{1}{(d-1)!}]$-algebra, or $R$ has property {\bf E}.
\end{enumerate} 
 Let $P$ be a projective $R$-module of rank $d$ with trivial determinant. Fix an isomorphism $\chi:R\iso \wedge^{d} P$.
The Euler class $e_{d}(P,\chi)$ is well-defined.
\et

\proof
Let $\beta:P\sur J'$ be another surjection such that $J'$ is an ideal of $R$ of height $d$.

By \cite[Lemma 3.0]{brs3}, there exists an ideal $I\subset R[T]$ of height $d$ and a surjection $\phi(T):P[T]\sur I$ such that $I(0)=J$,
$\phi(0)=\ga$ and $I(1)=J'$,  $\phi(1)=\gb$. 

Let $N=(I\cap R)^2$. Then $\hh(N)\geq d-1$ and therefore, $\dd(R/N)\leq 1$. By Serre's splitting theorem, $P/NP$ is a free $R/N$-module of rank $d$. On the other hand, by the same reasoning, $P[T]/IP[T]$ is a free $R[T]/I$-module of rank $d$.

We can choose an isomorphism $\tau:(R/N)^{d}\iso P/NP$ such that $\wedge^{d}\tau=\chi\otimes R/N$.
This choice of $\tau$ gives us a basis of $P/NP$, which in turn induces  a basis of the free module 
$P[T]/IP[T]$. Using this basis of $P[T]/IP[T]$ and the surjection $\phi(T):P[T]\sur I$, we obtain a surjection $\omega:(R[T]/I)^{d}\sur I/I^2$. Note that, due to the choice of the basis, $\omega(0)=\omega_J:(R/J)^{d}\sur J/J^2$, and $\omega(1)=\omega_{J'}:(R/J')^{d}\sur J'/J'^2$.

Now using Theorem \ref{chi}, we obtain an ideal $K$ of height $d$ with $K+J=R=K+J'$, and a 
surjection $\omega_K:(R/K)^{d}\sur K/K^2$ such that
$$(J,\omega_J)+(K,\omega_K)=(J',\omega_{J'})+(K,\omega_K) \text{ in } E^{d}(R).$$
This proves that $e_{d}(P,\chi)$ is well-defined.
\qed

\smallskip

\rmk
Note that, in \cite{brs3} it has been proved that the Euler class $e(P,\chi)$ is well-defined if $R$ is a commutative Noetherian ring  containing $\bq$. But essentially they require that $(d-1)!$ is invertible in $R$. We used that condition above. On the other hand, the Euler class is also well-defined  when $R$ is a smooth affine algebra over an infinite perfect field (the setup of \cite{brs1}).
We have stretched that to the condition \emph{``$R$ has property {\bf E}"} which will allow a bigger class of rings.

\medskip

Now that the Euler class is proved to be well-defined, the following theorem can easily be established using the arguments in \cite{brs1,brs3} verbatim.

\bt
Let $R,P$ be as in the above theorem. Then, $P\simeq Q\op R$ for some $R$-module $Q$ if and only if $e_d(P,\chi)=0$ in $E^d(R)$.  
\et

\rmk
Let us now comment on property {\bf E} and the Bass-Quillen conjecture. As mentioned above, if $R$ is a regular ring of dimension $2$, then by  Murthy's result \cite{mu1}, $R$ has property {\bf E}. Let $R$ be a regular ring of dimension $d\geq 3$. Thanks to the works of   Lindel \cite{li} and Popescu \cite{po}, 
we know that $R$ has property {\bf E} in the following cases: (a) $R$ contains a field, (b) The local rings of $R$ are either unramified or excellent Henselian.

In particular, as indicated in \cite{po}, the regular rings in the following theorem have property 
{\bf E} and therefore we can apply the Euler class theory to such rings. 

\bt\label{new1}
Let $R$ be a regular ring of dimension $d\geq 3$ which is either of the following: 
\begin{enumerate}
\item
$R$ is an $A$-algebra for some Dedekind domain $A$ such that for every $\mathfrak{p}\in \Spec(A)$
with $\mathfrak{p}R\neq R$, the ring $R/\mathfrak{p}R$ is regular and the quotient field of $A/\mathfrak{p}$ is perfect;
\item
In particular,  for every prime integer $p$, either $p$ is
a unit in $R$, or $R/pR$ is regular;
\item  
$R$ is smooth and of finite type over $\bz_p$ (the $p$-adic integers, $p$ a prime).
\end{enumerate}
Let $P$ be a projective $R$-module of rank $d$ with trivial determinant. Fix an isomorphism $\chi:R\iso \wedge^{d} P$. Then, $P\simeq Q\op R$ for some $R$-module $Q$ if and only if $e_d(P,\chi)=0$ in $E^d(R)$.
\et

\rmk
If $d=2$, then  the above theorem also holds without any regularity assumption.

\smallskip

In  Section \ref{cancel}, we shall  give a much simpler criterion for a projective $R$ module of rank $d$ to have a 
free summand, where $R$ is a $d$-dimensional affine $\bz$-algebra.

\section{The weak Euler class group}\label{wecg}
In this   section we recall the definitions of the weak Euler class groups from \cite{brs2,brs3}.
The notation $\mu(-)$ stands for the minimal number of generators.

\definition\label{defweak1}
({\bf The weak Euler class group  $E_0^d(R)$: $R$ smooth affine}) Let $R$ be a smooth affine domain of dimension $d\geq 2$
over an infinite perfect field $k$. Let $G_1$ be the free abelian group on the set $B_1$ of  maximal ideals of $R$.  Let $J=m_{1}\cap \cdots \cap m_r$, where $m_i$ are distinct  maximal
ideals of $R$.  We associate
$(J):= \sum_{1}^{r}(m_i)\in G_1$. Let $S_1$ be
the set of elements $(J)$ of $G_1$ such that $\mu(J)=d$.    Let $H_1$ be the 
subgroup of $G_1$ generated by $S_1$ .  The weak
Euler class group $E_0^d(R)$ is defined as $E_0^d(R):=G_1/H_1$.

\medskip

\definition\label{defweak2}
({\bf The weak Euler class group  $\wt{E}_0^d(R)$: $R$ Noetherian})  Let $R$ be a commutative Noetherian ring of dimension $d\geq 2$.
Let $G_2$ be the free abelian group on the set of ideals 
${\mathcal N}$ where ${\mathcal N}$ is an 
${\mathcal M}$-primary ideal for some maximal ideal ${\mathcal M}$ of height
$d$ such that  $\mu({\mathcal N}/{\mathcal N}^2)=d$ 
Now let $J\subset R$ be an ideal of height $d$ such that $\mu(J/J^2)=d$.
Let 
$J=\cap_{i}{\mathcal N}_{i}$ be  the (irredundant) primary decomposition
of $J$. We associate  $(J):=\sum_{i} ({\mathcal N}_{i})$ (in  $G_2$)
Let $H_2$ be the subgroup of $G_2$ generated by the set of elements $(J)$,
where $J$ is an ideal of height $d$  and $\mu(J)=d$.
The weak Euler class group of $R$  is 
$\wt{E}_0^d(R)\stackrel{\rm def}{=}G_2/H_2$.

\medskip

The same proof as in the case of the Euler class group will yield:

\bp
Let $R$ be a smooth affine domain of dimension $d\geq 2$
over an infinite perfect field $k$. Then $E_0^d(R)\simeq \wt{E}_0^d(R)$.
\ep

From now on, the weak Euler class group will be denoted as $E^d_0(R)$.
There is an obvious canonical surjective morphism $\psi:E^d(R)\sur E_0^d(R)$ taking an Euler cycle
$(J,\omega_J)$ to $(J)$. 

\bd
({\bf The weak  Euler class of a projective $R$-module:}) Let $R$ be a commutative Noetherian ring of dimension $d\geq 2$.
Let 
$P$ be a projective $R$-module of rank $d$ such that
$ R \simeq \wedge^d (P)$. Let 
$\varphi: P \sur J$ be a surjection where $J$ is an ideal of height $d$. The 
weak  Euler class of $P$ is defined as $e_d(P):=(J)\in E_0^d(R)$. 
\ed

It is easy to see that $e_d(P)$ is  
well-defined whenever
the Euler class $e_d(P,\chi)$ is so, for any fixed orientation $\chi:R\iso \wedge^d(P)$. Therefore,
for rings $R$ as in \ref{ecp}, $e_d(P)$ is well-defined. For such rings, one can easily prove the 
relevant set of results as in \cite{brs2,brs3}. We record one of them below.

\bt\label{weak}
Let $R$ be a ring as in \ref{ecp} and $P$ be a projective $R$-module of rank $d$ with trivial determinant. Assume that $P$ is stably free (or, more generally, $[P]=[Q\op R]$ for some $R$-module $Q$ of rank $d-1$). Then, $e_d(P)=0$.
\et

\proof
This is essentially contained in \cite{brs2,brs3}. Still, let us give a sketch. Assume that $P$ is stably isomorphic to $Q\op R$. By \cite[6.7]{brs3}, there exists an ideal $J$ of height $\geq d$
and surjections $\alpha:P\sur J$, $\gb:Q\op R\sur J$. If $J=R$, then $P$ has a free summand and $e_d(P,\chi)=0$ in $E^d(R)$ for any $\chi:R\iso \wedge^d(P)$. If $J$ is proper, then, as $Q\op R$ has a free summand, it follows from \cite[Theorem 1]{mo} that $\mu(J)=d$ and therefore, $(J)=0$ in $E_0^d(R)$. As $e_d(P)=(J)$, the result follows.
\qed

\section{Some exact  sequences}\label{etm}
We first emphasize that, unless mentioned otherwise, in this section we take $R$ to be  commutative and
Noetherian, \emph{without any further assumption}. 

Let $\dd(R)=d\geq 2$. The orbit space $Um_{d+1}(R)/E_{d+1}(R)$ has a group structure, thanks to the works of Vaserstein ($d=2$) \cite{suva} and van der Kallen ($d\geq 2$) \cite{vdk}. The orbit space $Um_{d+1}(R)/SL_{d+1}(R)$ also has a group structure and this group is a natural quotient of $Um_{d+1}(R)/E_{d+1}(R)$. In this short section we recall the definitions of maps from these groups to the Euler class groups.

\smallskip

\notation
Let $(J,\omega_J)\in E^d(R)$ and let $\ol{u}\in (R/J)^{\ast}$. Choose any $\sigma\in GL_d(R/J)$ with determinant $\ol{u}$. The notation $(J,\ol{u}\omega_J)$ stands for 
the Euler cycle $(J,\omega_J\sigma)$.  Note that this is well-defined as a local orientation
of $J$ is an equivalence class under the action of $SL_d(R/J)$.  

\smallskip

\subsection{Dimension two:} Let $\dd(R)=2$. Let $[a_1,a_2,a_3]\in Um_3(R)/SL_3(R)$.
Using elementary transformations we may assume that the height of the ideal  $J=(a_1,a_2)$ is $2$.
Let $\omega_J$ be the global orientation induced by $(a_1,a_2)$.
 Let  $\theta:R^3\sur R$ be given by $e_i\mapsto a_i$, $1\leq i\leq 3$, where $e_1,e_2,e_3$ are the standard basis vectors of $R^3$. Let $P=\text{ker}(\theta)$. Then we have a natural orientation $\chi:R\iso \wedge^2(P)$ (see \cite[p. 214]{brs3}) and a straightforward computation (\emph{op. cit.}) yields  that $e(P,\chi)=(J,\ol{a_3}\omega_J)$. The association $\varphi:[a_1,a_2,a_3]\mapsto e_d(P,\chi)=(J,\ol{a_3}\omega_J)$ is 
proved to be a morphism of groups in \cite{brs3} and the following sequence is exact:
\begin{equation}\label{1}
1\lra Um_3(R)/SL_3(R)\stackrel{\varphi}{\lra} E^2(R)\lra E_0^2(R)\lra 0
\end{equation}
Similarly, we have a morphism $\phi':Um_3(R)/E_3(R)\lra E^2(R)$ and an exact sequence:
\begin{equation}\label{2}
Um_3(R)/E_3(R)\stackrel{\phi}{\lra} E^2(R)\lra E_0^2(R)\lra 0
\end{equation}

\rmk
The maps $\phi,\phi'$ are well-defined because the Euler class $e_d(P,\chi)$ is so. Note that, as 
$d=2$, we do not need any additional assumption on $R$ for the Euler class to be well-defined.

\subsection{Higher dimensions:} Let $\dd(R)\geq 3$. Let $[a_1,\cdots,a_{d+1}]\in Um_{d+1}(R)/E_{d+1}(R)$.
As before, we may assume that the ideal $J=(a_1,\cdots,a_d)$ has height $d$. Let $\omega_J$ be the global orientation of $J$ induced by $a_1,\cdots,a_d$. It has been proved in \cite{dz} that the association $\phi:[a_1,\cdots,a_{d+1}]\mapsto (J,\ol{a_{d+1}}\omega_J)$ is well-defined, and is a morphism. Further, the following sequence of groups is exact:
\begin{equation}\label{3}
Um_{d+1}(R)/E_{d+1}(R)\stackrel{\phi}{\lra} E^d(R)\lra E_0^d(R)\lra 0
\end{equation}

\section{Some relation with cancellation}\label{cancel}
As mentioned above, in this section we intend to give a very straightforward criterion for a projective $R$ module of rank $d$ to have a 
free summand, where $R$ is a $d$-dimensional affine $\bz$-algebra (without any smoothness assumption). This criterion is a result of the cancellative nature of the free module $R^d$ (see \cite{suva}). We shall provide some finer analysis of this phenomenon.  Moreover, we think it is better to give a more general and unified treatment, in the form of following set of results.

In what follows we make the following assumption: 

\noindent
($*$) \,\,\,\emph{Let $A$ be a commutative Noetherian ring such that for any finite type $A$-algebra $S$ of dimension 
$\leq 2$, the free module $S^2$ is cancellative.}

\bt
Let $R$ be a ring of $\text{dim}(R)=d\geq 2$, which is of finite type over $A$, where $A$ as in ($*$). 
Then, the canonical morphism $\psi:E^d(R)\lra E^d_0(R)$ is an isomorphism.
\et

\proof
We first assume that $\dd(R)=d=2$. In this case, by hypothesis, $R^2$ is cancellative and therefore,
$Um_{d+1}/SL_{d+1}(R)$ is trivial. The exact sequence \ref{1} from Section \ref{etm} gives us the desired result.

We now assume that $d\geq 3$.
As $\psi$ is already surjective, we need to check its injectivity. Let $(J,\omega_J)\in E^d(R)$ be
such that $\psi(J,\omega_J)=(J)=0$ in $E_0^d(R)$. Then, from the exact sequence \ref{3} in Section \ref{etm}, it follows that there is a unimodular row $[a_1,\cdots,a_{d+1}]$ such that:
\begin{enumerate}
\item
$K:=(a_1,\cdots,a_d)$ has height $d$;
\item
$(K,\ol{a_{d+1}}\omega_K)=(J,\omega_J)$, where $\omega_K$ is the global orientation induced by $a_1,\cdots,a_{d}$.
\end{enumerate}
Therefore, it is enough to prove that $(K,\ol{a_{d+1}}\omega_K)=0$. Note that we can also ensure that
$\hh(a_1,\cdots,a_{d-2})=d-2$. Write $S=R/ (a_1,\cdots,a_{d-2})$. 

Now $\ol{a_{d+1}}\omega_K$ is induced by $K=(a_{1},\cdots,a_{d+1}a_d)+K^2$. Let tilde denote reduction modulo $ (a_1,\cdots,a_{d-2})$. Then, in $S$ we have $\wt{K}=(\wt{a}_{d-1},\wt{a}_{d})$ and
$$\wt{K}=(\wt{a}_{d-1},\wt{a}_{d+1}\wt{a}_d)+\wt{K}^2.$$

If $\dd(S)=2$, by hypothesis we have $S^2$ is cancellative and arguing as in the first paragraph of this proof, we obtain that $\wt{K}=(\wt{b},\wt{c})$ such that $\wt{b}-\wt{a}_{d-1}\in \wt{K}^2$, and
$\wt{c}-\wt{a}_{d+1}\wt{a}_d\in \wt{K}^2$. We can then find suitable preimages $b,c$ in $K$ so that  $K=(a_1,\cdots,a_{d-2},b,c)$ with $b-a_{d-1}\in K^2$ and $c-a_{d+1}a_d\in K^2$. This proves that 
$(K,\ol{a_{d+1}}\omega_K)=0$. The case when $\dd(S)\leq 1$ follows from standard general position arguments.
\qed

\smallskip

For the following result, we cannot invoke the Euler class of a projective module.
However, the underlying ideas are from the theory of Euler classes (as in \cite[3.3, 3.4]{brs3}).

\bt\label{new2}
Let $R$ be a ring of dimension $d\geq 2$, which is of finite type over $A$, where $A$ as in ($*$). Let $P$ be a projective $R$-module of rank $d$ with $\wedge^d(P)\iso R$. Assume that there is an 
$R$-linear surjection $\alpha:P\sur J$ such that $J$ is generated by $d$ elements. Then 
$P\simeq Q\oplus R$ for some $R$-module $Q$.
\et

\proof
Let $\chi: R\iso \wedge^d(P)$ be an isomorphism. Observe that $P/JP$ is a free $R/J$-module of rank $d$.
We choose an isomorphism $\sigma:(R/J)^d\iso P/JP$ such that $\wedge^d(\sigma)=\chi\otimes_R R/J$. Write $\overline{\alpha}=\alpha\otimes_R R/J$. Let $\omega_J:(R/J)^d\sur J/J^2$ be the composite:
$$(R/J)^d\stackrel{\sigma}{\iso} P/JP\stackrel{\overline{\alpha}}{\sur} J/J^2$$
Let $\omega_J$ correspond to $J=(b_1,\cdots,b_d)+J^2$. On the other hand, it is given that $J$ is generated by $d$ elements. This implies that $(J)=0$ in $E^d_0(R)$, and by the above theorem, 
$(J,\omega_J)=0$ in $E^d(R)$. 

\smallskip

\noindent
\emph{Case 1.} Let $d=2$. In this case, the Euler class of a projective module is well-defined
by Theorem \ref{ecp} and we have $e_2(P,\chi)=(J,\omega_J)=0$. Therefore, $P\simeq Q\oplus R$ for some 
$R$-module $Q$.

\smallskip

\noindent
\emph{Case 2.} Let $d\geq 3$. We have, $(J,\omega_J)=0$ in $E^d(R)$. 
Therefore, there exist
$c_1,\cdots,c_d\in J$ such that  $J=(c_{1},\cdots,c_{d})$ and $b_{i}\equiv c_{i}$ mod $J^2$ for $i=1,\cdots,d$. 
As the height of $J$ is $d$, we can find $\lambda_1,\cdots,\lambda_{d-1}\in R$ such that 
 $\text{dim}(R/(c'_1,\cdots,c'_{d-1})\leq 1$,
where $c'_i=c_i+\lambda_i c_d$ for $i=1,\cdots,d-1$. Set $b'_i=b_i+\lambda_i b_d$ for $i=1,\cdots,d-1$. Therefore, we have, $J=(c'_{1},\cdots,c'_{d})$ and $b'_{i}\equiv c'_{i}$ mod $J^2$ for $i=1,\cdots,d$. Observe that the operations just performed correspond to an elementary matrix. We can then alter $\sigma$ by this elementary matrix to obtain $\sigma'$ so that the composite $\omega'_J$
$$(R/J)^d\stackrel{\sigma'}{\iso} P/JP\stackrel{\overline{\alpha}}{\sur} J/J^2$$ gives rise to 
$J=(b'_1,\cdots,b'_d)+J^2$. The upshot of this series of arguments is that, without loss of generality, we may assume that $\text{dim}(R/(c_1,\cdots,c_{d-1})\leq 1$ to start with. We assume this and proceed.

Consider the polynomial algebra $R[T]$ and the following ideal in $R[T]$:
$$I=(c_1,\cdots,c_{d-1},T+c_d).$$
As $\text{dim}(R[T]/I)=\dd(R/(c_1,\cdots,c_{d-1})\leq 1$, the projective $R[T]/I$-module $P[T]/IP[T]$ is free of rank $d$.
We choose an isomorphism $\theta:P[T]/IP[T]\iso (R[T]/I)^d$ such that $\wedge^d \theta=\chi(T)^{-1}\otimes_{R[T]} R[T]/I$. Substituting $T=0$ we observe that $\theta(0)^{-1}$ and 
$\sigma$ differ by an automorphism $\overline{\delta}\in SL_{d}(R/J)$. But 
$SL_d(R/J)=E_d(R/J)$ and we can lift $\overline{\delta}$ to $\delta\in E_d(R)$. 
Taking $\theta':=\delta(T)\theta:P[T]/IP[T]\iso  (R[T]/I)^d$ we ensure that $\theta'(0)^{-1}=\sigma$.

Let $\beta:R[T]^d\sur I$ be the surjection which sends $e_i$ to $c_i$ for $i=1,\cdots,d-1$, and $e_d$ to $T+c_d$. We then have the induced surjection 
$$\gamma:=\overline{\beta}\theta:P[T]/IP[T]\sur I/I^2.$$
Observe that $\gamma(0)=\overline{\alpha}$. We can now apply a result of  Mandal \cite[Theorem 2.1]{m2} to obtain a surjection
$\eta:P[T]\sur I$. Substituting $T=1-c_d$, we are done.
\qed

\bc\label{new3}
Let $R$ be a ring of dimension $d\geq 2$ and $P$ be a projective $R$-module of rank $d$ with trivial determinant. Assume that there is an $R$-linear surjection $\alpha:P\sur J$ such that
$J$ is an ideal of height $d$ which is generated by $d$ elements. Then $P\iso Q\oplus R$ for some $Q$ in the following cases:
\begin{enumerate}
\item
$R$ is of finite type over $\bz$.
\item
$R$ is of finite type over a field $k$ which satisfies one of the following:
\begin{enumerate}
\item
$p\neq \text{char}(k)$, $c.d._p(k)\leq 1$;
\item
$p=\text{char}(k)$ and $k$ is perfect;
\end{enumerate}
\end{enumerate}
\ec

\proof
(1) Take $A=\bz$ and note that $A$ satisfies the condition ($*$) due to  a result of Suslin-Vaserstein\cite[Corollary 18.1, Theorem 18.2]{suva}.

\smallskip

(2) Take $A=k$ as above and note that $A$ satisfies the condition ($*$) due to a result of Suslin 
\cite[Theorem 2.4]{su},
modified by Parimala \cite[Theorem 3.1]{rv}.
\qed

\smallskip

With the notations set up at the beginning of this section, one can easily  prove the following 
``addition" and ``subtraction" principles.

\bt
Let $R$ be a ring of dimension $d\geq 2$, which is of finite type over $A$ (where $A$ as in ($*$)). Let $I$ and $J$ be two 
comaximal ideals of $R$, each of height $d$. If two of $I$, $J$, $I\cap J$ are generated by $d$ elements, then so is the third.
\et

\section{The Euler class group $E^d(R[T])$}\label{ecgt}

We may recall that in \cite{d1}, the $d\,\textsuperscript{th}$ Euler class group $E^d(R[T])$ has been defined where $R$ is a commutative Noetherian ring of dimension $d\geq 3$ \emph{containing the field of rationals}. Recall also that for a Noetherian ring $A$ of dimension $\gd$, and an integer 
$n\geq \frac{1}{2}(\gd+3)$, the $n$-th Euler class group $E^n(A)$ had already been defined in \cite{brs4}. Therefore, if one takes
a Noetherian ring $R$ of dimension $d\geq 4$, the theory of $E^d(R[T])$ is available from \cite{brs4} by taking $A=R[T]$, $\gd=d+1$, and $n=d$. This was precisely done in \cite{my}. However, this approach has a couple of disadvantages. First of all, it leaves out the case $d=3$. Secondly, it does not allow one to define the Euler class of a projective $R[T]$-module $P$ of rank $d$ (together with $\chi:R[T]\iso \wedge^d(P)$).
 In order to define the Euler class of $(P,\chi)$, one has to define an Euler cycle $(I,\omega_I)\in E^d(R[T])$ with the \emph{local orientation} $\omega_I$ as an equivalence class of surjections induced by the action of
$SL_d(R[T]/I)$ on all surjections from $(R[T]/I)^d\sur I/I^2$. We did so in \cite{d1}. In comparison, the definition in \cite{brs4} uses the
action of the elementary group $E_d(R[T]/I)$ (which may be a proper subgroup of $SL_d(R[T]/I)$).

In concurrence with \cite{d1}, we shall define $E^d(R[T])$ in this section where $R$ is as described in the remark below.

\rmk\label{hypo}
In what follows, our hypothesis for the ring $R$ will be as follows.
We consider $R$ to be a regular ring of dimension $d\geq 3$ with either of the following assumptions:
\begin{enumerate}
\item
$R$ is smooth and of finite type over a finite field $k$; or  
\item
$R$ contains an infinite field $k$.
\end{enumerate}
Note that $R$ has property {\bf E} (see Definition \ref{prope}).

\medskip

\rmk
Let $I\subset R[T]$ be an ideal of height $d$. 
If $R$ contains an infinite field $k$, then by \cite[Lemma 3.3]{brs1}, there exists $\lambda\in k$ such that
$I(\lambda)=R$ or $\hh(I(\lambda))=d$. Changing $T$ by $T-\lambda$, without loss of generality, one can assume that $I(0)=R$ or $\hh(I(0))=d$. We used it heavily in \cite{d1}. Here, when $R$ is of finite type over a finite field, we have to use
some different arguments than \cite{d1}.
We illustrate one such key situation below for the convenience of the reader.

\bt\label{sl}
Let $R$ be as in \ref{hypo}. Let $I\subset R[T]$ be an ideal of height $d$ with $\mu(I/I^2)=d$. Let $\ga,\gb$ be two surjections from $(R[T]/I)^d$ to $I/I^2$ be such that there exists 
$\sigma\in SL_d(R[T]/I)$ with $\ga\sigma=\gb$. Assume that $\ga$ can be lifted to a surjection 
$\theta:R[T]^d\sur I$. Then, $\gb$ can also be lifted to a surjection 
$\phi:R[T]^d\sur I$.
\et

\proof
We first show that we can lift $\gb$ to a surjection $\gamma:R[T]^d\sur I/(I^2T)$.
We tackle  two cases separately.

\smallskip

\noindent
{\it Case 1.} We assume that $R$ is smooth and of finite type over a finite field $k$.

If $I(0)=R$, then by \cite[3.9]{brs1} we can lift $\gb$ to a surjection $\gamma:R[T]^d\sur I/(I^2T)$. We now show that we can do the same if $I(0)$ is proper. Note that $\hh(I(0))\geq d-1$ and therefore, $\dim(R/I(0))\leq 1$. We have 
$\alpha(0):(R/I(0))^d\sur I(0)/I(0)^2$, $\gb(0):(R/I(0))^d\sur I(0)/I(0)^2$ and $\ga(0)\sigma(0)=\gb(0)$. Now, 
by \cite[Remark after 16.2]{suva}, $K_1Sp(R/I(0))=SK_1(R/I(0))$ whereas, by \cite[Theorem 9.11]{swan} we have
$K_1Sp(R/I(0))=0$. Therefore, $SK_1(R/I(0))=0$. Since
by \cite[Theorem 9.1]{swan}, stable rank of $R/I(0)\leq 2$, and $d\geq 3$, it follows that $SL_d(R/I(0))=E_d(R/I(0))$.
As $E_d(R)\lra E_d(R/I(0))$ is surjective, we can find $\tau\in E_d(R)$ which is a lift of $\sigma(0)$. Clearly,
$\theta(0)\tau:R^d\sur I(0)$ is a  lift of $\gb(0)$. By \cite[3.9]{brs1}, $\gb$ can be lifted to a surjection $\gamma:R[T]^d\sur I/(I^2T)$.

\smallskip

\noindent
{\it Case 2.} We assume that $R$ contains an infinite field $k$. Using \cite[Lemma 3.3]{brs1}, if necessary, we can replace $T$ by $T-\lambda$ for a suitable $\lambda\in k$ so that $I(0)=R$ or $\hh(I(0))=d$. Then, $\dd(R/I(0))=0$ and therefore $SL_d(R/I(0))=E_d(R/I(0))$.
Rest is  same as the last part of Case 1.

\smallskip

From both the cases, we have $\gamma:R[T]^d\sur I/(I^2T)$ which is a lift of $\gb$.
We now move to $R(T)$ and note that $\dd(R(T)/IR(T))=0$. Therefore, we have $SL_d(R(T)/IR(T))=E_d(R(T)/IR(T))$ and the matrix $\sigma\ot R(T)$ has a lift to $\sigma'\in E_d(R(T))$. Then, $(\theta\ot R(T))\sigma':R(T)\sur IR(T)$ is a lift of $\beta\ot R(T)$. We can
now apply \cite[Proposition 4.9]{bk} to conclude that $\gb$ has a lift to a surjection $\phi:R[T]^d\sur I$.
\qed

Using \cite[Proposition 4.9]{bk} and the addition, subtraction principles from \cite{drs1}, it is easy to prove the following addition and subtraction principles. Basic line of arguments are the same as in \cite[Propositions 4.2, 4.3]{d1}.

\bt\label{add} (Addition principle) Let $R$ be as in \ref{hypo}. Let $I,J\subset R[T]$ be two comaximal ideals, each of height $d$.
Let $I=(f_1,\cdots,f_d)$ and $J=(g_1,\cdots,g_d)$. Then $I\cap J=(h_1,\cdots,h_d)$ where $h_i-f_i\in I^2$ and $h_i-g_i\in J^2$ for $i=1,\cdots,d$.
\et

\bt\label{sub} (Subtraction principle) Let $R$ be as in \ref{hypo}. Let $I,J\subset R[T]$ be two comaximal ideals, each of height $d$.
Let $I=(f_1,\cdots,f_d)$ and $I\cap J=(h_1,\cdots,h_d)$ where $h_i-f_i\in I^2$ for $i=1,\cdots,d$. Then $J=(g_1,\cdots,g_d)$  with $h_i-g_i\in J^2$ for $i=1,\cdots,d$.
\et

Let $I\subset R[T]$ be an ideal of height $d$ such that $\mu(I/I^2)=d$. The group $SL_d(R[T]/I)$ acts on the set of all surjections  $(R[T]/I)^d\sur I/I^2$. By a \emph{local orientation}  of $I$ we mean an equivalence class under that action. Theorem \ref{sl} tells us that if a surjection
$\ga:(R[T]/I)^d\sur I/I^2$ has a lift to surjection $R[T]^d\sur I$, then the same happens for  any 
$\gb$ in the equivalence class of $\ga$.
A local orientation $[\ga]$ is called \emph{global} if
$\ga$ has a lift to a surjection $R[T]^d\sur I$. 
As before, by abuse of notations, we shall identify a surjection $\ga:(R[T]/I)^d\sur I/I^2$ with its equivalence class.

\smallskip

\definition\label{def}
({\bf The Euler class group $E^d(R[T])$:}) Let $R$ be as in \ref{hypo}. Let $G$ be the free abelian group on the set $B$ of pairs 
$(\CI,\omega_{\CI})$
where $\CI\subset R[T]$ is an ideal of height $d$ such that 
Spec\,$(R[T]/\CI)$ is connected, $\CI/\CI^{2}$ is generated by $d$ elements and 
$\omega_{\CI}:(R[T]/\CI)^{d}\twoheadrightarrow \CI/\CI^{2}$ is a local orientation of 
$\CI$. 
Let $I\subset R[T]$ be an ideal of height $d$ and $\omega_{I}$ a local 
orientation of $I$. Now $I$ can be decomposed uniquely as 
$I=\CI_{1}\cap\cdots\cap \CI_{r}$, where the $\CI_{i}$'s are ideals of $R[T]$ 
of height $d$, pairwise comaximal and Spec\,$(R[T]/\CI_{i})$ is connected for 
each $i$. Clearly $\omega_{I}$ induces  local orientations
$\omega_{\CI_{i}}$ of $\CI_{i}$ for $1\leq i \leq r$. By 
$(I,\omega_{I})$ we mean the element $\Sigma(\CI_{i},\omega_{\CI_{i}})$ of $G$.
Let $H$ be the subgroup of $G$ generated by set $S$ of pairs $(I,\omega_{I})$,
where $I$ is an ideal of $R[T]$ of height $d$ generated by $d$ elements 
and $\omega_{I}$ is a global 
orientation of $I$ given by the set of generators of $I$. 
We define the $d\,\textsuperscript{th}$ Euler class group of $R[T]$, denoted
by $E^d(R[T])$, to be $G/H$.

We now state a lemma from \cite[4.1]{k}.

\begin{lemma}\label{mphil}
Let $G$ be a free abelian group with basis $B=(e_i)_{i\in {\mathcal I}}$ . Let $\sim$  be an equivalence relation on
$B$. Define $x\in G$  to be ``reduced" if $x=e_1+\cdots +e_r$ and $e_i\not=e_j$ for $i \not=j$.
%For $x\in G$ with $x=e_1+\cdots +e_r$ define``support"  of $x$ to be the set ${\rm supp}(x):=\{e_1,\cdots,e_r\}$. 
Define $x\in G$ to be ``nicely reduced" if $x=e_1+\cdots +e_r$ is such that  $e_i\not\sim e_j$ for $i \not=j$.
Let $S\subset G$ be such that 
\begin{enumerate}
\item
Every element of $S$ is nicely reduced.
\item
Let $x,y\in G$ be such that each of $x,y,x+y$ is nicely reduced. If two of $x,y,x+y$ are in $S$, then so is the third.
\item
Let $x\in G\setminus S$ be nicely reduced and let ${\mathcal J}\subset {\mathcal I}$ be  finite. Then there exists $y\in G$ with the following properties : (i) $y$ is nicely reduced; (ii) $x+y\in S$; (iii) $y+e_j$ is nicely reduced $\forall j\in{\mathcal J}$. 
\end{enumerate}

Let $H$ be the subgroup of $G$ generated by $S$. If $x \in H$ is nicely reduced, then $x \in  S$.
\end{lemma}

We are now ready to prove:

\bt\label{zero}
Let $R$ be as in \ref{hypo}.
Let $I$ be any ideal of $R[T]$ of height $d$ such that $\mu(I/I^2)=d$, and 
$\omega_I:(R[T]/I)^{d}\sur I/I^2$ be a local orientation.
Assume that the image of $(I,\omega_I)$ is trivial in $E^{d}(R[T])$. Then $\omega_I$ is global.
In other words, there is an
$R[T]$-linear  surjection $\theta : (R[T]/I)^{d}\sur I$ such that $\theta $ lifts $\omega_I$.
\et

\proof
We take $G$ to be the free abelian group generated by $B$, as defined in (\ref{def}). Define a relation $\sim$ on $B$ as: $(K,\omega_K)\sim (K',\omega_{K'})$ if 
$K=K'$. Then it is an equivalence relation. 

Let $S\subset G$ be as in (\ref{def}). In view of the above lemma, it is enough to show that the three
conditions in (\ref{mphil}) are satisfied. Condition (1) is clear, almost from the definition.
The addition and subtraction principles (\ref{add}, \ref{sub}) will yield  condition (2). Finally, applying the moving lemma 
\cite[Lemma 2.12]{d1}, it is clear that (3) is also satisfied.
\qed

\subsection{Homotopy invariance}
Let  $F$ be a functor from the category of commutative Noetherian rings to the category of groups.
Then, $F$ is said to enjoy \emph{homotopy invariance} if for a regular ring $R$, $F(R)\iso F(R[T])$.
Although the Euler class groups lack nice functorial properties, we can nevertheless ask whether the
homotopy invariance holds. It is easy to see that the obvious canonical morphism 
$\varphi: E^d(R)\lra E^d(R[T])$   is injective. The best result that we have is the following:

\bt\label{inv}
Let $R$ be a regular ring  of dimension $d\geq 3$ which is essentially of finite type over an infinite field $k$. Then the canonical morphism $\varphi: E^d(R)\lra E^d(R[T])$ is an isomorphism.
\et

\proof
The same method as in \cite[Proposition 5.7]{d1} works. Instead of \cite[Theorem 3.8]{brs1}, one has to use
\cite[Theorem 4.12]{d3}, as we did not assume $k$ to be perfect here.
\qed

\medskip

\quest
Let $R$ be smooth and of finite type over a \underline{finite} field $k$. Is the canonical morphism $\varphi: E^d(R)\lra E^d(R[T])$  an isomorphism?

\section{Two results of Asok-Fasel}\label{asofa}

Let $X=\Spec(R)$ be a smooth affine $k$-algebra of dimension $d\geq 2$ where $k$ is an infinite perfect field. As we have already collected the necessary results in the previous sections, here we shall just remark on the  changes required in the arguments of \cite[Section 3]{af}.

We first rewrite a  crucial result from \cite[Lemma 3.1.11]{af}. We stick to the notations in their article and we do not explain any notation or terminology. This lemma is about the consistency of the definition of a map $\wt{\psi}_d:Q_{2d}(R)\l E^d(R)$
where $R$ is a smooth affine algebra of dimension $d\geq 3$ over an infinite perfect field $k$. The reader will notice that we confine ourselves to the `top' case $d=n$ here, as the final results are about that setup. Also, we include the case $d=2$ below which is not treated in \cite[Lemma 3.1.11]{af}.

\bl
Let $X=\Spec(R)$ be a smooth affine $k$-scheme of dimension $d\geq 2$, where $k$ is  an infinite perfect field.
Then map $\wt{\psi}_d:Q_{2d}(R)\l E^d(R)$ is well-defined.
\el

\proof
There are two cases.

\smallskip

\noindent
\emph{Case 1.} Assume that $d\geq 3$.
Follow the proof of \cite[Lemma 3.1.11]{af} and at the end of the proof use Theorem \ref{inv}.

\smallskip

\noindent
\emph{Case 2.} Assume that $d=2$. (Following their notations) we have $\omega_0,\omega_1$
who are in the same homotopy orbit as $v$. Now, $\omega_i$ defines $(N_{\mu_{i}},\omega_{N_{\mu_{i}}})
\in E^2(R)$ for $i=0,1$.  

Following the proof, we finally have 
an ideal $I\subset R[T]$ of height $2$ and a surjection $\omega_I:(R[T]/I)^2\sur I/I^2$ such that 
$I(0)=N_{\mu_{0}}$, $I(1)=N_{\mu_{1}}$ and $\omega_I$ restricted to $T=i$ gives  $\omega_{N_{\mu_{i}}}$
for $i=0,1$. By \cite[4.3, 4.6]{brs1}, the elements  $(N_{\mu_{0}},\omega_{N_{\mu_{0}}})$ and
$(N_{\mu_{1}},\omega_{N_{\mu_{1}}})$ are equal in $E^2(R)$.
\qed

\smallskip

\rmk
Ideas similar to the proof of Case 2 were used in \cite[Theorem 3.10]{dtz1}.

\smallskip

Using the above lemma, one can extend  \cite[Proposition 3.1.12]{af} and eventually obtain the following theorem (cf. \cite[Theorem 3.1.13]{af}).

\bt
Let $X=\Spec(R)$ be a smooth affine $k$-scheme of dimension $d\geq 2$, where $k$ is  an infinite perfect field. The Segre class morphism $s: E^d(R)\l [X,Q_{2d}]_{{\mathbb A}^1}$ is an isomorphism.
\et

As a consequence,  the following improvement of \cite[Theorem 3.2.1]{af} is obtained, the proof being the same as in \cite{af}. Also note that, since we have included the case $d=2$ in the above results, one
does not need to treat this case separately as has been done in \cite{af}.

\bt
Let $k$ be an infinite perfect field of characteristic unequal to $2$. Let $X=\Spec(R)$ be a smooth affine $k$-scheme  of dimension $d\geq 2$. Then the canonical morphism $E^d(R)\lra \wt{CH}^d(X)$ is an isomorphism.
\et

Now we come to the comparison between the weak Euler class group $E_0^d(R)$ and the Chow group (of zero cycles) $CH^d(X)$. 

\bt
Let $k$ be an infinite perfect field and Let $X=\Spec(R)$ be a smooth affine $k$-scheme  of dimension $d\geq 2$. Then the canonical morphism $E_0^d(R)\lra CH^d(X)$ is an isomorphism in the following cases:
\begin{enumerate}
\item
$d=2$;
\item
$d\geq 3$ and $\text{Char}(k)\neq 2$.
\end{enumerate}
\et

\proof
If we make the blanket assumption that $\text{Char}(k)\neq 2$, then the same proof as in \cite{af} works as we can use the exact sequence (2) from Section \ref{etm} for the case $d=2$ as well.

However, for $d=2$ one can give a direct proof without invoking the Chow-Witt group (where the characteristic assumption is required). Quite a long time ago,  Bhatwadekar had told me that it is not  difficult to prove that $E_0^2(R)\iso CH^2(X)$.
Somehow we never discussed his ideas. Let me try to give some arguments here.

Each of the following morphisms is surjective:
$$ E^2_0(R)\sur CH^2(X)\sur F^2K_0(R)$$
Let $(J)\in E^2_0(R)$ be such that $J$ is reduced of height $d$ and its image in $CH^2(X)$ is trivial.
  We have $[R/J]=0$ in $F^2K_0(R)$. Now, by \cite[2.2]{mu2}, there is a projective $R$-module of rank $2$ and a surjection $\alpha:P\sur J$ such that 
$$[P]-[R^2]=-[R/J] \text{ in } F^2K_0(R).$$
As $[R/J]=0$, it follows that $P$ is stably free. But we have the weak Euler class of $P$, 
$e_0(P)=(J)\in  E^2_0(R)$. By Theorem \ref{weak}, it follows that $(J)=0$ in $E^2_0(R)$.
Therefore, $E^2_0(R)\iso CH^2(X)\iso F^2K_0(R)$.
\qed

\end{document}